\documentclass[11pt,english]{article}
\usepackage[T1]{fontenc}
\usepackage[latin9]{inputenc}
\usepackage{amsmath}
\usepackage{graphicx,pstricks}
\usepackage{amssymb}
\usepackage[arrow, matrix, curve]{xy}
\usepackage[english]{babel}

\usepackage[footnotesize]{caption} 

\makeatletter

\usepackage{epsfig}\usepackage{graphics}

\newtheorem{thm}{Theorem}

\def\0{{\bf 0}}

\def\phi{\varphi}

\def\T{\T}

\newcommand{\bk}[1]{\ensuremath{\langle#1\rangle}}
\newcommand{\fine}{~\hfill~$\square$\newline}
\newcommand{\e}{\varepsilon}

\catcode`@=11 \@addtoreset{equation}{section} \catcode`@=12

\begin{document}


\title{Time to absorption  for a heterogeneous neutral competition model}

\author{Claudio Borile\thanks{\noindent Dipartimento di Fisica ``G. Galilei'', Universit\`a di
Padova, via Marzolo 8, I-35151 Padova, Italy,\newline\noindent \texttt{borile@pd.infn.it}},\,  
Paolo Dai Pra\thanks{\noindent Dipartimento di Matematica,
Universit\`a di Padova, via Trieste 63, I-35121 Padova, Italy,\newline\noindent \texttt{daipra@math.unipd.it}},\,
Markus Fischer\thanks{\noindent Dipartimento di Matematica,
Universit\`a di Padova, via Trieste 63, I-35121 Padova, Italy,\newline \texttt{fischer@math.unipd.it}},\, 
Marco Formentin\thanks{\noindent Dipartimento di Fisica ``G. Galilei'', Universit\`a di
Padova, via Marzolo 8, I-35151 Padova, Italy,\newline \texttt{marco.formentin@rub.de}}, and 
Amos Maritan\thanks{\noindent Dipartimento di Fisica ``G. Galilei'', Universit\`a di
Padova, via Marzolo 8, I-35151 Padova, Italy,\newline \texttt{maritan@pd.infn.it}}
}

\maketitle

\begin{abstract} 
Neutral models aspire to explain biodiversity patterns in ecosystems where species difference can be neglected, as it might occur at a specific trophic level, and perfect symmetry is assumed between species. Voter-like models capture the essential ingredients of the neutral hypothesis and represent a paradigm for other disciplines like social studies and chemical reactions. In a system where each individual can interact with all the other members of the community, the typical time to reach an absorbing state with a single species scales linearly with the community size. Here we show, by using a rigorous approach within a large deviation principle and confirming previous approximate and numerical results, that in a heterogeneous voter model the typical time to reach an absorbing state scales exponentially with the system size, suggestive of an asymptotic active phase.
\end{abstract}

\bigskip 

\noindent\textit{\textbf{keywords:}}{ Voter Model with disorder, Neutral models of biodiversity, Large deviations, Stochastic dynamics with quenched disorder}

\section{Introduction} \label{sect:intro}


Models of interacting degrees of freedom are nowadays widely spread in different scientific disciplines---from Physics and Mathematics to Biology, Ecology, Finance and Social Sciences---, and more than ever in the last few years there has been a growing effort in connecting the phenomenology observed at a macroscopic level with a simplified ``microscopic'' modeling of very disparate complex systems. Clearly, this idea is extremely appealing to statistical physicists and can provide a good benchmark for developing new ideas and methods. A famous and particularly successful example of this approach, which reconciles interdisciplinarity and pure research in statistical physics, can be found in the ecological literature in the so-called \textit{neutral theory of species diversity}, that aims at giving a first null individual-based modeling of the dynamic competition among individuals of different species in the same trophic level of an ecosystem \cite{Hubbell,Volkov03,Azaele06,vallade2003analytical,alonso2006merits}. The neutral hypothesis finds its mathematical equivalent in the voter model (VM) \cite{Liggett} and its generalizations \cite{AlHammal05}, which, in turn, is equivalent to the well-known Moran model in genetics \cite{Kimura}. This model has been deeply studied and has gradually become a paradigmatic example of non-equilibrium lattice models. It is conceptually simple but nevertheless has a very rich phenomenology with applications in many different scientific areas \cite{Henkel,Dornic01,Durrett94,Castellano09b,Blythe07}. Despite the fact that the original formulation of the VM can be exactly solved in any spatial dimension \cite{Liggett}---fact that contributed greatly to its rise---, any slight modification made in order to improve the realism of the model complicates drastically its analysis.

Among the possible modifications of the original VM, there has been a recent interest in studying the asymptotic behavior of the VM in the presence of quenched random-field-like disorder, whose motivations span from ecology 
\cite{Pigolotti10,borile2013effect,tilman1994habitat} to social modeling \cite{Masuda10,Masuda11} through models of chemical reactants \cite{Krapivsky} and more fundamental research \cite{Odor}. A particularly interesting problem is to assess the typical time needed by a finite-size system to reach one of the absorbing states of the model; depending on the particular interpretation of the model, that would mean the typical time for the extinction of a species in an ecosystem, or the reaching or not of a consensus on a particular topic in a society. In all these cases, it is known that heterogeneities, in the habitat of an ecosystem or in the ideologies of groups of people, play a major role in shaping the global dynamics of the complex system.  It has been shown \cite{Masuda10} that a quenched (random-field-like) disorder creating an intrinsic preference of each individual for a particular state/opinion hinders the formation of consensus, hence favoring coexistence. In the context of neutral ecology, this corresponds to a version of the VM in which at each
location there is an intrinsic preference for one particular species, leading to mixed states lasting for times that grow exponentially with system size \cite{Pigolotti10,Masuda10}.

Here, we propose a rigorous mathematical development of a disordered VM intended as a general model of neutral competition in a heterogeneous environment. Supporting the previous findings \cite{Pigolotti10,Masuda10} based on computational investigations or approximation arguments, we will show that a heterogenous environment indeed favors significantly the maintenance of the active state, and the typical time needed to reach an absorbing phase passes from a power-law dependence in the system's size, typical of the neutral theories, to an exponentially long time, signature of an asymptotic active phase. This will be achieved by setting up a large deviation principle for the considered model and will thus provide a first attempt of an extreme value theory for systems with multiple symmetric absorbing states.\\

\section{Macroscopic limit} The state of the system is described by a vector of spins $\eta = (\eta_1,\eta_2,\ldots,\eta_N) \in \{0,1\}^N$. The  random environment consists of $N$ independent and identically distributed random variables $h_1,h_2,\ldots,h_N$,  taking the values $0$ and $1$ with probability, respectively, $1-q$ and $q \in (0,1)$. Moreover, let $\rho \in [0,1]$ be a given parameter. While the random environment remains constant, the state $\eta$ evolves in time according to the following rules:
\begin{itemize}
\item each site $1,2,\ldots,N$ has its own independent random clock. A given site $i$ after a waiting time with exponential distribution of mean $1$ chooses at random, with uniform probability, a site $j \in \{1,2,\ldots,N\}$.\\
\item If $\eta_j = h_i$, then the site $i$ updates its spin from $\eta_i$ to $\eta_j$. If $\eta_j \neq  h_i$, then the site $i$ updates its spin from $\eta_i$ to $\eta_j$ with probability $\rho$, while it keeps its spin $\eta_i$ with probability $1-\rho$.
\end{itemize}
Thus, the site $i$ has a preference to agree with sites whose spins equal its local field $h_i$. For $\rho = 1$, this effect is removed, and we obtain the standard Voter model. Note that, by symmetry, there is no loss of generality in assuming $q \geq 1/2$, as we will from now on.

In more formal terms, for every realization of the random environment, the spins evolve as a continuous-time Markov chain with generator $L_N$ acting on a function $f:\{0,1\}^N\rightarrow \mathbb{R}$ according to
\begin{align}\label{spin}
	L_Nf(\eta) := \sum_{i=1}^N\frac{1}{N}&\sum_{j=1}^N\mathbb{I}_{h_i=\eta_j}\left(f(\eta^{j\rightarrow i})-f(\eta)\right)\nonumber\\
&+\rho\mathbb{I}_{h_i\neq\eta_j}\left(f(\eta^{j\rightarrow i})-f(\eta)\right),
\end{align}
where $\mathbb{I}_A$ denotes the indicator function of the set $A$ and $\eta^{j\rightarrow i}$ is the configuration obtained from
$\eta$ by replacing the value of the spin at the site $i$ with that of the spin at the site $j$. This Markov chain has two absorbing states, corresponding to all spin values equal to zero and all equal to one. We denote by $T_N$ the random time needed to reach one of the two absorbing states. \\

It is useful to review the main properties of the model in the case $\rho = 1$. 
In this case, the dynamics are independent of $q$ and  the unique order parameter for the model is given by
$K_N :=  \sum_{i=1}^N \eta_i$,
i.e., the number of spins with value $1$. It is easy to check, using the generator \eqref{spin}, that $K_N$ evolves as a random walk on $\{0,1,\ldots,N\}$: if $K_N = k$, then it moves to either $k+1$ or $k-1$ with the same rate $\frac{(N-k)k}{N}$. By standard arguments on birth and death processes (see e.g. \cite{lambert2008}), one shows that $\bk{T_N} \sim N \ln(2)$ as $N \rightarrow + \infty$: the mean absorption time grows linearly in $N$.

Consider now the general case $\rho \leq 1$. Here the system is described in terms of two integer-valued order parameters, namely $\displaystyle{\sum_{i=1}^Nh_i\eta_i}$ and $\displaystyle{\sum_{i=1}^N(1-h_i)\eta_i}$, that will be convenient to properly scale as follows:
\begin{equation*}
\begin{split}
m_N^+ & := m_N^+(\eta) := \frac{1}{N}\sum_{i=1}^Nh_i\eta_i \\
m_N^- & := m_N^-(\eta) :=\frac{1}{N}\sum_{i=1}^N(1-h_i)\eta_i
\end{split}
\end{equation*}
Note that the pair $(m_N^+, m_N^-)$ belongs to the subset of the plane
$\{(x,y) \in [0,1]^2 : x+y \leq 1\}$.
Note, however, that when the limit as $N \rightarrow +\infty$ is considered, $m^+_N \leq  \frac{1}{N}\sum_{i=1}^Nh_i \rightarrow q$, where this last convergence follows from the law of large numbers. Similarly, $m^-_N \leq  \frac{1}{N}\sum_{i=1}^N(1-h_i) \rightarrow 1- q$. Thus, limit points of the sequence $(m_N^+, m_N^-)$ belong to $[0,q] \times [0,1-q]$.
Given an initial, possibly random, state $\eta(0)$ for the dynamics of $N$ spins, we denote by $m_N^{\pm}(t)$ the (random) value at time $t$ of the order parameters $m_N^{\pm}$. In what follows we also denote by $\mu_N$ the distribution of $\eta(0)$.\medskip\\

\begin{thm} \label{limdyn}
Assume there exists a non-random pair $(\bar{m}^+, \bar{m}^-) \in [0,q] \times [0,1-q]$ such that, for every $\epsilon>0$,
\[
\lim_{N \rightarrow +\infty} \mu_N\left( \left| m_N^{\pm}(0) - \bar{m}^{\pm} \right| >  \epsilon \right) = 0.
\]
Then the stochastic process $(m^+(t), m^-(t))_{t \geq 0}$ converges in distribution to the unique solution of the following system of ODEs:
{\renewcommand*{\arraystretch}{1.3}
\begin{equation}\label{eq}
\left\{\begin{array}{lll}
\dot{m}^+ &= &-\rho m^+(1-m^--m^+)\\
 & &+(q-m^+)(m^++m^-)\\
\dot{m}^- & = &-m^-(1-m^--m^+)\\
& &+\rho(1-q-m^-)(m^++m^-) \\ m^{\pm}(0) & = & \bar{m}^{\pm}
\end{array}\right.
\end{equation}
}
\end{thm}

\noindent{\textbf{Proof:}}
Denote by $\mathcal{G}$ the generator of the semigroup associated to the deterministic evolution \eqref{eq}, i.e.,
\[
\mathcal{G} f(m^+,m^-) := V^+(m^+,m^-) \frac{\partial f}{\partial m^+}   +  V^-(m^+,m^-) \frac{\partial f}{\partial m^-},
\]
with
\[
\begin{split}
V^+(m^+,m^-) = &  -\rho m^+(1-m^--m^+) \\ & +(q-m^+)(m^++m^-) \\
V^-(m^+,m^-) = &  -m^-(1-m^--m^+) \\ & + \rho(q-m^-)(m^++m^-)
\end{split}
\]
Let $f: [0,1]^2 \rightarrow \mathbb{R}$. By direct computation one finds that
\[
L_N [f(m^+_N, m^-_N)](\eta)
\]
depends on $\eta$ only through $m^+_N,m^-_N$, which implies that the process $(m_N^+(t), m_N^-(t))_{t \geq 0}$ is a Markov process, whose associated semigroup has a generator $\mathcal{G}_N$ that can be identified by the identity
\[
L_N [f(m^+_N, m^-_N)](\eta) = [\mathcal{G}_N f] (m^+_N(\eta), m^-_N(\eta)),
\]
which yields
\begin{equation} \label{gn}
\begin{split}
	&\mathcal{G}_{N}f(x,y)\\
	&:= N \Bigl( \left(q-x\right) (x+y) \left(f(x+\tfrac{1}{N},y) - f(x,y)\right)) \\
	&\quad+ \rho x \left(1-(x+y)\right)  \left(f(x-\tfrac{1}{N},y) - f(x,y)\right)  \\
	&\quad+ \rho  \left(1-q-y\right)  (x+y) \left(f(x,y+\tfrac{1}{N}) - f(x,y)\right)  \\
	&\quad+ y  \left(1-(x+y)\right)  \left(f(x,y-\tfrac{1}{N}) - f(x,y)\right)  \Bigr) .
\end{split}
\end{equation}

Moreover, if $f$ is smooth with bounded derivatives, one checks that
\[
\lim_{N \rightarrow +\infty}\sup_{(m^+,m^-) \in [0,1]^2} \left| \mathcal{G}_N f(m^+,m^-) -  \mathcal{G} f(m^+,m^-)\right| = 0.
\]
The conclusion then follows by a standard result of convergence of Markov processes, cf. \cite{ethier2009markov}, Ch. 4, Corollary 8.7.
\fine

{This first theorem formalizes and extends a useful result for the infinite size system, already obtained by means of different techniques in some previous works \cite{Masuda10,borile2013effect}. It {is a dynamic law of large numbers that} quantifies 
the deterministic evolution of the order parameters as obtained
from the limiting dynamics described by $L_N$ neglecting fluctuations. The stability analysis of the fixed points of Eq.\eqref{eq} provides some immediate results on the global dynamics of the model in the infinite size limit:}
For $\rho = 1$, equations \eqref{eq} trivialize: the only relevant variable is $m  = m^+ + m^-$, which satisfies $\dot{m} = 0$. This is simply the macroscopic consequence of the fact that $K_N = N (m^+_N+ m_N^-)$ evolves as a symmetric random walk. The picture changes as $\rho <1$. When $\rho<1$, the system \eqref{eq} has three equilibrium points:
\begin{enumerate}
\item $(m^+, m^-)=(q,1-q)$, which represents a limiting behavior where all the spins equal 1;\\
\item $(m^+, m^-)=(0,0)$, which is the case with all spins equal to 0;\\
\item $(m^+, m^-)= \left( \frac{q(1+\rho)-\rho}{(1+\rho)(1-\rho)}, \rho \frac{q(1+\rho)-\rho}{(1+\rho)(1-\rho)} \right)$.
\end{enumerate}
It is easily checked that equilibrium 3 lies inside $[0,q] \times [0,1-q]$, hence is admissible, 
if and only if the condition
\begin{equation} \label{qcon}
\rho < \frac{1-q}{q}
\end{equation}
holds (remember we are assuming $q \geq 1/2$). The stability analysis of the three equilibria is also easily done: for $\frac{1-q}{q} < \rho < 1$ equilibrium 1 is stable, and attracts all initial conditions except $(0,0)$, which is an unstable equilibrium, while for $\rho < \frac{1-q}{q}$ both $(q,1-q)$ and $(0,0)$ are unstable, and the stable equilibrium 3 emerges, attracting all initial conditions except the unstable equilibria. Note that, for $q = 1/2$, only this second regime exists.
Thus, in the case $q > 1/2$ and $\frac{1-q}{q} < \rho < 1$, the asymmetric disorder stabilizes the equilibrium $(q,1-q)$; lower values of $\rho$ increase the effects of the disorder, so that a new stable equilibrium appears.\\

\section{Large deviations and time to absorption} {In order to get information on the behavior of the system when the total number of ``individuals'' $N$ is large but finite, we need to go beyond the law of large numbers in Theorem \ref{limdyn}. In particular, our next aim is } to show that, whenever equilibrium 3 is present for the macroscopic dynamics \eqref{eq}, the absorption time for the microscopic system grows exponentially in $N$. To this end, we use the Freidlin and Wentzell theory for randomly perturbed dynamical systems (see \cite{freidlin2012random}). This theory, based on finite time Large Deviations, yields asymptotic estimates characterizing the long-time behavior of the perturbed system (here the microscopic system described by $(m^{+}_{N},m^{-}_{N})$) as the noise intensity tends to zero (equivalent here to $N\to \infty$). See also \cite{den2008large,touchette2009large} for an introduction to Large deviations.\\
For simplicity, we assume $q = 1/2$, so that equilibrium 3 exists for every $\rho <1$. For ${\bf x} = (x,y) \in [0,1/2]^2$ set
\[
\begin{split}
l_1({\bf x}) & =  \rho x (1-x-y)  \\  r_1({\bf x }) & = (1/2 - x)(x+y) \\
l_2({\bf x}) & = y (1-x-y)  \\  r_2({\bf x }) & = \rho (1/2 - y)(x+y)
\end{split}
\]
Notice that the vector field $b({\bf x}) = (b_1({\bf x}), b_2({\bf x}))$ defined by $b_i({\bf x}) := r_i({\bf x }) - l_i({\bf x })$ appears in \eqref{eq}, the equation of the macroscopic dynamics, which we interpret as the unperturbed dynamical system. Define the family of point measures, parametrized by ${\bf x} \in [0,1/2]^2$:
\begin{equation*}
	\mu_{\boldsymbol{x}} := r_{1}(\boldsymbol{x})\delta_{(1,0)} + l_{1}(\boldsymbol{x})\delta_{(-1,0)} + r_{2}(\boldsymbol{x})\delta_{(0,1)} + l_{2}(\boldsymbol{x})\delta_{(0,-1)},
\end{equation*}
where $\delta$ indicates Dirac measure. Then the generator $\mathcal{G}_{N}$ in \eqref{gn} can be rewritten in a diffusion-like form as
\begin{multline*}
	\mathcal{G}_{N}(f)(\boldsymbol{x}) = N \int_{\mathbb{R}^{2}\setminus\{0\}} \left( f(\boldsymbol{x} + \tfrac{1}{N}\boldsymbol{\gamma}) - f(\boldsymbol{x})\right) \mu_{\boldsymbol{x}}(d\boldsymbol{\gamma}) \\
	= \left\langle b(\boldsymbol{x}), \nabla f(\boldsymbol{x})\right\rangle \\  +  N \int \left( f(\boldsymbol{x} + \tfrac{1}{N}\boldsymbol{\gamma}) - f(\boldsymbol{x}) - \frac{1}{N} \left\langle \boldsymbol{\gamma}, \nabla f(\boldsymbol{x})\right\rangle \right) \mu_{\boldsymbol{x}}(d\boldsymbol{\gamma}),
\end{multline*}
where $\langle \cdot , \cdot \rangle$ denotes the scalar product in $\mathbb{R}^2$. Let $H\!: \mathbb{R}^{2}\times \mathbb{R}^{2} \rightarrow \mathbb{R}$ be the Hamiltonian associated with the operators $\mathcal{G}_{N}$, $N\in \mathbb{N}$:
\[
	H(\boldsymbol{x},\boldsymbol{\alpha}) :=  \left\langle b(\boldsymbol{x}), \boldsymbol{\alpha}\right\rangle + \int \left( \exp\left(\langle \boldsymbol{\gamma}, \boldsymbol{\alpha} \rangle\right) - 1 - \left\langle \boldsymbol{\gamma}, \boldsymbol{\alpha}\right\rangle \right) \mu_{\boldsymbol{x}}(d\boldsymbol{\gamma}).
\]
It follows that
\[
	H(\boldsymbol{x},\boldsymbol{\alpha}) = \sum_{i=1}^{2} \left[r_{i}(\boldsymbol{x})\left(e^{\alpha_{i}} - 1 \right) + l_{i}(\boldsymbol{x})\left(e^{-\alpha_{i}}- 1 \right) \right].
\]
Let $L$ be the Legendre transform of $H$, given by
\[
	L(\boldsymbol{x},\boldsymbol{\beta}) := \sup_{\boldsymbol{\alpha}\in \mathbb{R}^{2}}\left\{\left\langle \boldsymbol{\beta}, \boldsymbol{\alpha}\right\rangle - H(\boldsymbol{x},\boldsymbol{\alpha})\right\}.
\]
It is easy to show that
\begin{equation} \label{EqLagrangian}
	L(\boldsymbol{x},\boldsymbol{\beta}) = \tilde{L}(l_{1}\bigl(\boldsymbol{x}),r_{1}(\boldsymbol{x});\beta_{1}\bigr) + \tilde{L}\bigl(l_{2}(\boldsymbol{x}),r_{2}(\boldsymbol{x});\beta_{2}\bigr),
\end{equation}
where $\tilde{L}: [0,\infty)^{2}\times \mathbb{R} \rightarrow [0,\infty]$ is given by
\[
\begin{split}
	\tilde{L}(l,r;\beta) & =  \sup_{\alpha\in \mathbb{R}} \left\{ \beta\cdot\alpha - r\cdot (e^{\alpha}-1) - l\cdot (e^{-\alpha}-1)\right\} \\
	& = \beta\log\left(\tfrac{\beta + \sqrt{\beta^{2}+4rl}}{2r}\right) - \sqrt{\beta^{2}+4rl} + l + r,
\end{split}
\]
taking appropriate limits for the boundary cases $l=0$ or $r = 0$. In particular, $\tilde{L}(l,r;\beta) = \infty$ if and only if either $l=0$ and $\beta < 0$ or $r=0$ and $\beta > 0$.
The {\it Lagrangian} $L$ in \eqref{EqLagrangian} allows to define the {\it action functional}: for $T>0$, $\phi:[0,T] \rightarrow \mathbb{R}^2$, set
\begin{equation} \label{ExActionFnct}
	S_{T}(\phi) :=
		\int_{0}^{T} L\bigl(\phi(t),\dot{\phi}(t)\bigr)dt,
\end{equation}
where $S_{T}(\phi)$ is meant to be equal to $+\infty$ if $\phi$ is not absolutely continuous. The action functional controls the {\it quenched} Large Deviations of the stochastic process $(m_N^+(t), m_N^-(t))_{t \geq 0}$: if $B_{\phi}$ is a small neighborhood of a trajectory $\phi:[0,T] \rightarrow \mathbb{R}^2$, $h = (h_1,h_2,\ldots,h_N)$ is a realization of the random environment, and $P_h$ is the law of the Markov process generated by \eqref{spin} for $h$ {\it fixed}, then for almost every realization $h$
\[
\frac{1}{N} \log P_h\left[ (m^+(t), m^-(t))_{t \in [0,T]} \in B_{\phi} \right] \simeq - S_{T}(\phi)
\]
for $N$ large. This fact falls within the range of the Freidlin-Wentzell Large Deviations results (see \cite{freidlin2012random}), although several modifications of the original proof are needed here, following \cite{Budhiraja11}.

As shown in \cite{freidlin2012random}, the control of the Large Deviations provides control on the hitting times of subsets of the state space $[0,1/2]^2$ of the process $(m_N^+(t), m_N^-(t))_{t \geq 0}$, in particular of the time $T_N$ needed to reach the absorbing states. Denote by ${\bf z}$ the stable equilibrium for the macroscopic dynamics:
\[
{\bf z} = \left(\frac{1}{2(1+\rho)}, \frac{\rho}{2(1+\rho)} \right).
\]
For ${\bf x} \in [0,1/2]^2$, define the {\it quasi-potential} by
\[
V({\bf x}) := \inf \{ S_T(\phi) : T > 0, \phi(0) = {\bf z}, \phi(T) = {\bf x} \}.
\]
Let $D$ be a domain in $[0,1/2]^{2}$ containing ${\bf z}$ with smooth boundary $\partial D$ such that $\partial D \subseteq (0,1/2)^{2}$ and the vector field $b({\bf x})$ is directed strictly inside $D$. Let $\tau_{N}$  denote the first time the process $(m^{+}_{N},m^{-}_{N})$  hits the complement of $D$. By construction, $\tau_{N} \leq T_{N}$.\medskip\\

\begin{thm} \label{quasipotential}
For every ${\bf x} \neq {\bf z}$ we have $V({\bf x}) >0$. Moreover, for almost every realization of the environment $h$, every $\e>0$,
\[
\lim_{N \rightarrow +\infty} P_h \left(e^{N(V_{\partial D}-\e)} \leq \tau_N \leq e^{N(V_{\partial D}+\e)} \right) = 1
\]
where
\[
V_{\partial D} := \min\left\{ V({\bf x}) : {\bf x}\in \partial D \right\} > 0.
\]
\end{thm}
\noindent{\textbf{Proof:} }
In order to show that $V({\bf x}) >0$ for every ${\bf x} \neq {\bf z}$, it suffices to check that, for every $\delta_{0} > 0$ small enough, $\inf_{\bf x\in \partial B_{\delta_{0}}(\boldsymbol{z})} V(\bf x) > 0$.

Set ${\bf r_{\ast}}:=(\frac{\rho}{4(1+\rho)},\frac{\rho}{4(1+\rho)})$; thus ${\bf r_{\ast}} = (r_{i}({\bf z}),l_{i}({\bf z}))$, $i\in \{1,2\}$. Let $l,r > 0$. Then $\tilde{L}(l,r;\beta)$ as a function of $\beta\in \mathbb{R}$ is smooth, non-negative, strictly convex with minimum value zero attained at $\beta = r-l$ and of super-linear growth. Second order Taylor expansion around $\beta = r -l$ yields
\[
	\tilde{L}(l,r;\beta) = \frac{1}{2(r+l)}\left(\beta - (r-l)\right)^{2} + \mathcal{O}\left(\left(\beta - (r-l)\right)^{3}\right).
\]
It follows that for every $\delta_{\ast} > 0$ small enough there are a constant $c > 0$ and a continuous function $\underline{L}\!: \overline{B_{\delta_{\ast}}(\bf{r}_{\ast})}\times \mathbb{R} \rightarrow [0,\infty)$ such that $\tilde{L}(l,r;\beta) \geq \underline{L}(l,r;\beta)$, $\underline{L}(l,r;.)$ is strictly convex with super-linear growth and for every $(l,r)\in \overline{B_{\delta_{\ast}}(\boldsymbol{r}_{\ast})}$,
\[
	\underline{L}(l,r;\beta) = c \left(\beta - (r-l)\right)^{2} \text{ if } \beta\in [-4\delta_{\ast},4\delta_{\ast}].
\]
Choose such $\delta_{\ast}$, $c$, $\underline{L}$. By continuity of the functions $r_{1}$, $l_{1}$, $r_{2}$, $l_{2}$, we can choose $\delta_{0} > 0$ such that $(l_{1}(\boldsymbol{x}),r_{1}(\boldsymbol{x})), (l_{2}(\boldsymbol{x}), r_{2}(\boldsymbol{x})) \in \overline{B_{\delta_{\ast}}(\boldsymbol{r}_{\ast})}$ for all $\boldsymbol{x}\in \overline{B_{\delta_{0}}(\boldsymbol{z}))}$. Recall that $b_{i} = r_{i} - l_{i}$. It follows that
\[
\begin{split}
	\inf_{\boldsymbol{x}\in \partial B_{\delta_{0}}(\boldsymbol{z})} &\bar{V}(\boldsymbol{x}) \\
	&\geq \inf \sum_{i=1}^{2}\int_{0}^{T} \underline{L}\left(l_{i}(\phi(t)),r_{i}(\phi(t)); \dot{\phi}_{i}(t)\right)dt,
\end{split}
\]
where the infimum on the right-hand side is over all $\phi\in \mathbf{C}_{a}([0,\infty), \overline{B_{\delta_{\ast}}(\boldsymbol{r}_{\ast})})$, $T > 0$ such that $\phi(0)=\boldsymbol{z}$, $\phi(T) \in \partial B_{\delta_{0}}(\boldsymbol{z})$. Using a time transformation argument analogous to that of Lemma~4.3.1 in \cite{freidlin2012random} and the convexity and super-linear  growth of $\underline{L}(l,r,\beta)$ in $\beta$, one finds that the infimum can be restricted to $\phi\in \mathbf{C}_{a}([0,\infty), \overline{B_{\delta_{\ast}}(\boldsymbol{r}_{\ast})})$ such that $|\dot{\phi}_{i}(t)|\leq 4\delta_{\ast}$ and $|\dot{\phi}(t)| = |b(\phi(t))|$ for almost all $t\in \mathbb{R}$. Thus
\[
	\inf_{\boldsymbol{x}\in \partial B_{\delta_{0}}(\boldsymbol{z})} \bar{V}(\boldsymbol{x}) \geq \inf \int_{0}^{T} c\cdot \bigl|b(\phi(t)) - \dot{\phi}(t)\bigr|^{2} dt.
\]
The Jacobian of $b$ at $\boldsymbol{z}$ has two strictly negative eigenvalues. Choosing, if necessary, a smaller $\delta_{\ast}$ and corresponding $c > 0$, $\underline{L}$, it follows that
\[
	\inf_{\boldsymbol{x}\in \partial B_{\delta_{0}}(\boldsymbol{z})} \bar{V}(\boldsymbol{x})
	\geq \inf
	\int_{0}^{T} c\cdot \bigl|Db(\boldsymbol{z})\phi(t) - \dot{\phi}(t)\bigr|^{2} dt > 0,
\]
which establishes the strict positivity of $V$ away from ${\bf z}$.

The second part of the assertion is established in a way analogous to the proofs of Theorems 4.4.1 and 4.4.2 in \cite{freidlin2012random}, see Section 5.4 therein.
\fine

Theorem~\ref{quasipotential} implies in particular that the time to reach any small neighborhood of the absorbing states grows exponentially in $N$ { for any $\rho<1$. This is a generalization of the Kramers's formula for the noise activated escape from a potential well \cite{hanggi1986escape}. This exponential behavior in $N$ suggests the existence of an active phase where both spin states / species, $0$ and $1$, coexist in the stationary state in the infinite size limit, $N\rightarrow \infty$.} \\

\section{Normal fluctuation} As seen in the previous sections, on a time scale of order $1$ the process $(m_N^+(t), m_N^-(t))_{t \geq 0}$ remains close to its thermodynamic limit: i.e., Eq.\eqref{eq}. In this section we consider the normal fluctuations around this limit. Suppose the assumptions of Theorem \ref{limdyn} are satisfied; moreover, for the sake of simplicity, we assume $q = 1/2$, and $(m^+,m^-) = {\bf z}$ with ${\bf z} = \left(\frac{1}{2(1+\rho)}, \frac{\rho}{2(1+\rho)} \right)$, so that the limiting dynamics starts in equilibrium. We define the fluctuation processes
\[
\begin{split}
x_N(t) & := \sqrt{N} \left(m^+_N(t) - m^+ \right) \\
y_N(t) & := \sqrt{N} \left(m^-_N(t) - m^- \right) .
\end{split}
\]
\medskip\\
\begin{thm}\label{fluctuations}
The stochastic process $(x_N(t), y_N(t))$ converges in distribution to a Gauss-Markov process $(X,Y)$ which solves the stochastic differential equation
{\renewcommand*{\arraystretch}{1.8}
\begin{equation}\label{diff}\left\{\begin{array}{ll}
d X=&\left(-\frac{1+\rho^2}{2(1+\rho)}X+\frac{\rho}{1+\rho}Y+\frac{1}{2}\mathcal{H}\right)dt\\
&+\frac{1}{\sqrt{2}}\sqrt{\frac{\rho}{1+\rho}}d B_1\\
d Y=&\left(-\frac{1+\rho^2}{2(1+\rho)}Y+\frac{\rho}{1+\rho}X-\rho\frac{1}{2}\mathcal{H}\right)dt\\
&+\frac{1}{\sqrt{2}}\sqrt{\frac{\rho}{1+\rho}}d B_2
\end{array}\right.
\end{equation}
}
Here, $B_i$, $i=1,2$ are two independent standard Brownian motions and $\mathcal{H}$ is a zero average standard Gaussian random variable, independent of $B_1,B_2$.
\end{thm}
\medskip
The proof of Theorem \ref{fluctuations} uses the method of convergence of generators as that of Theorem \ref{limdyn}, and is omitted. Unlike in Theorem \ref{limdyn}, the environment does not fully self-averages since
$\mathcal{H}$ is not identically equal to zero. The quenched random variable $\mathcal{H}$ in Theorem \ref{fluctuations} is due to the normal fluctuations of the environment $(h_1,h_2,\ldots,h_N)$.\\

\section{Discussion and conclusions}
It is well known that habitat heterogeneity impacts on biodiversity \cite{mcclain2010habitat}. At large scale, e.g. at regional or larger level, geomorphological changes may induce genetic isolation whereas at smaller scales the complexity induced by, for example, vegetation, sediment types, moisture and temperature leads to the coexistence of several species and to the emergence of niches. {To our knowledge, however, quantitative estimates of the relation between the degrees of heterogeneity and biodiversity and the time of coexistence of species have not been obtained. Here we have rigorously proved that even a small habitat disorder in a neutral competition-like model dramatically enhances the typical time biodiversity persists;
more specifically, we have shown that the typical time to loss of biodiversity, $\tau_N$, scales exponentially with the population size $N$, leading, for large size systems, to an unobservable long time scale beyond which extinction occurs. This is in contrast to what happens in absence of habitat heterogeneity, where the typical time to loss of biodiversity is typically small, growing as the system's size, $N$. We have also obtained the scaling exponent of $\tau_N$ in terms of a suitable \textit{quasi-potential} $V(\bf{x})$, that encodes the minimum ``cost'' of a trajectory to reach a given point $x$ of the phase-space. The consequences of these findings could be particularly relevant, for example, in conservation ecology: In a given area different species at the same trophic level compete for space and nutrients in a neutral fashion; for example, think of a tropical forest, where the neutral theory provides a very good null model \cite{Volkov03}.\\
Lastly, we have shown that the fluctuations around the metastable symmetric fixed point obey a Brownian motion dynamics with drift where the environmental disorder does not show self-averaging.}

\medskip


\textbf{\textit {Acknowledgments:}}
AM acknowledges Cariparo foundation for financial support. We thank Miguel Munoz  for useful discussions,  comments and suggestions.

\bibliographystyle{unsrt}


\end{document}